\documentclass[11pt]{article} 
\usepackage{latexsym}
\usepackage{amsfonts}
\usepackage{amsthm}
\usepackage{amsmath}
\usepackage{amssymb}

\usepackage{eucal}
\usepackage{graphicx}
\usepackage[all, poly, knot]{xy}

\input xy
\xyoption{all}
\xyoption{knot}
\xyoption{arc}

\input epsf.tex
\newdimen\epsfxsize
\newdimen\epsfysize
\def\qed{\vrule height5pt width3pt depth.5pt}

\newtheorem{thm}{Theorem}[section]
\newtheorem{cor}[thm]{Corollary}
\newtheorem{lem}[thm]{Lemma}

\theoremstyle{definition}

\newtheorem{rem}{Remark}[section]
\title{Vassiliev Invariants from Parity Mappings}

\author{H. A. Dye}

\begin{document}

\maketitle

\begin{abstract} Parity mappings from the chords of a Gauss diagram to the integers is defined. The parity of the chords is used to construct families of invariants of Gauss diagrams and virtual knots. One family consists of degree $n$ Vassiliev invariants. 
\end{abstract}

\section{Introduction}
\begin{figure}[[htb] \epsfysize = 1 in
\centerline{\epsffile{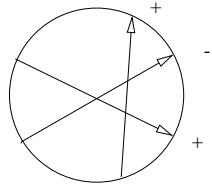}}
\caption{Gauss diagram for the trefoil}
\label{fig:gaussdiagram}
\end{figure}

A \textit{Gauss diagram} is a clockwise oriented circle with with $n$ signed, oriented chords. The orientation of each chord is indicated by an arrowhead and the sign of each chord is indicated with a $+$ or $-$.  A Gauss diagram is shown in figure \ref{fig:gaussdiagram}. The endpoint  marked by the arrowhead is referred to as the head and the other endpoint is referred to as the foot.

\begin{figure}[[htb] \epsfysize = 3 in
\centerline{\epsffile{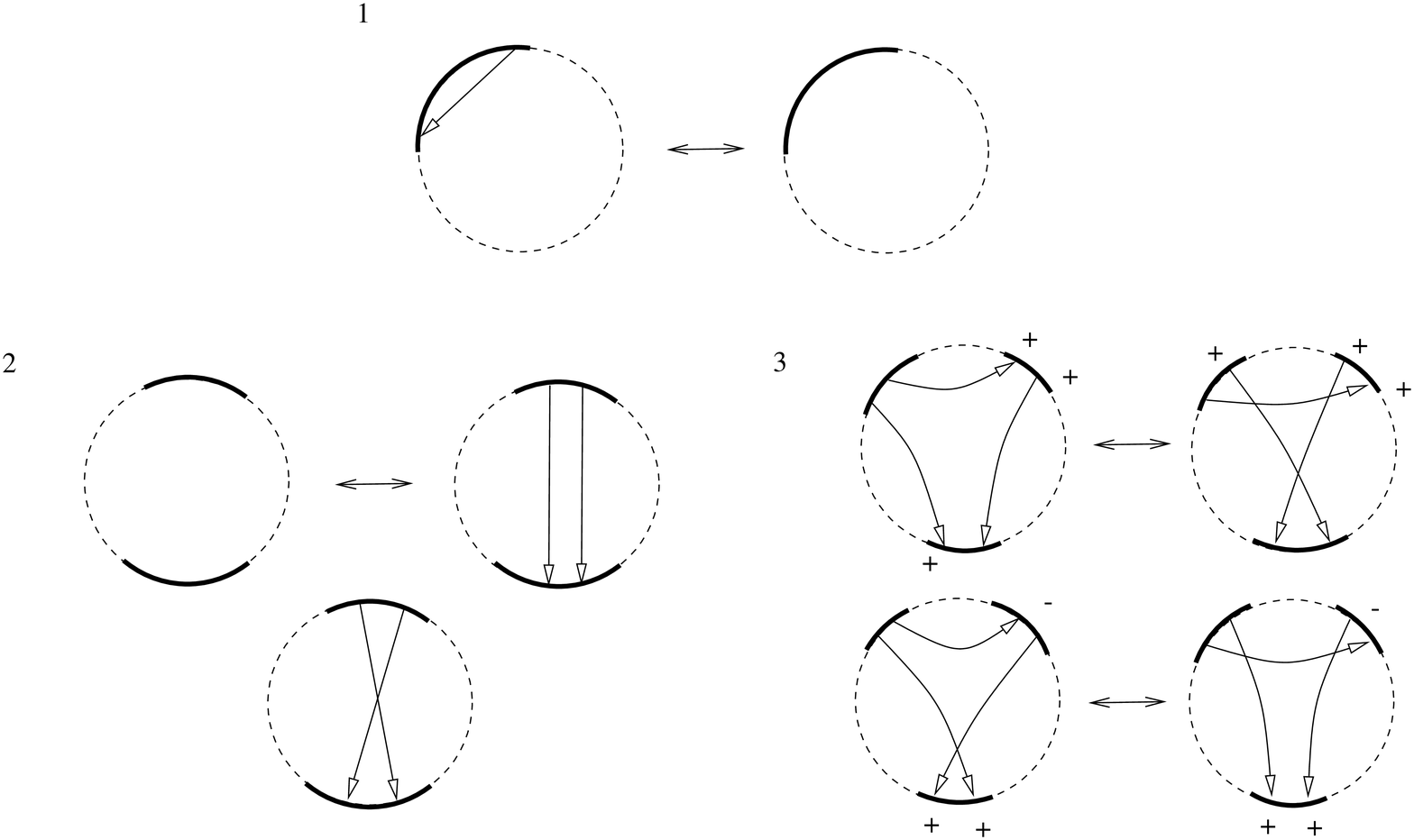}}
\caption{Moves on Gauss diagrams}
\label{fig:gaussmoves}
\end{figure}

We define a set of  three diagrammatic chord moves on Gauss diagrams. 
The \textit{single chord} move introduces a single chord that does not intersect any chords.
The \textit{two chord} move introduces two chords with opposite signs. These chords are positioned so that the heads are adjacent and the feet are adjacent. 
The \textit{triangle} move acts on three existing chords in a Gauss diagram. The
endpoints of the chords occur in pairs and the chords form a triangle. 
There are two versions of the triangle move and the version depends on the intersection of the chords (figure \ref{fig:gaussmoves}).
The (3,0) triangle move is an equivalence between a local configuration of chords with three intersections and a local configuration with zero intersections. The (2,1) triangle move is an equivalence between a local configuration of chords with two intersections and a local configuration of chords with one intersection.
The signs and orientation of the chords must be selected correctly. In the
 the triangle formed by the chords, one vertex consists of two heads, one vertex consists of two feet, and one contains a head and a foot. 
In figure \ref{fig:gaussmoves}, we show a correct selection of signs for the two versions of the triangle move. 
From these diagrams, we can modify the diagrams to obtain the other correct arrangements of orientation and sign. We change the sign and orientation of a pair of chords with adjacent heads (respectively feet) in both diagrams. 
The set of chord moves determine equivalence classes of Gauss diagrams. Two diagrams are equivalent if they are related by finite sequence of  chord moves.

\begin{figure}[[htb] \epsfysize = 2.5 in
\centerline{\epsffile{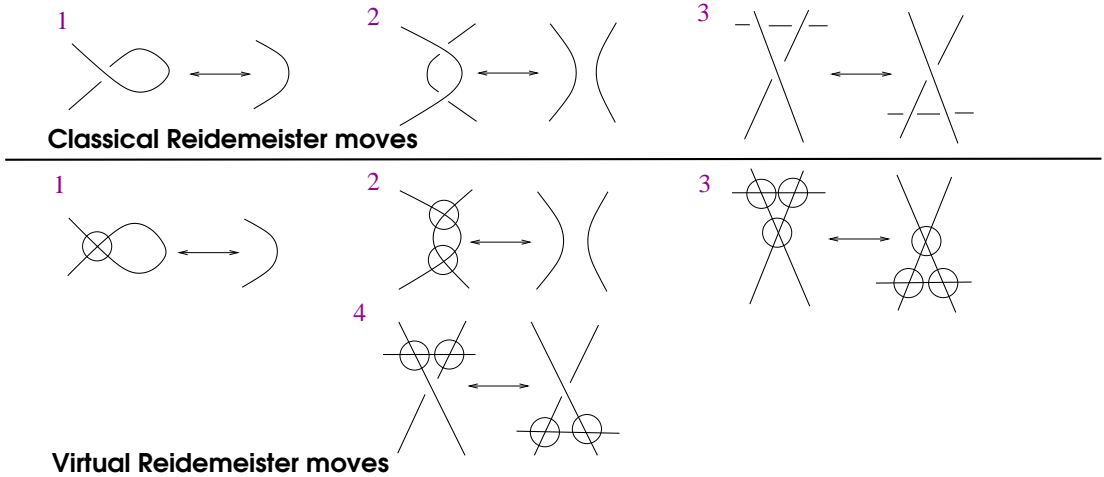}}
\caption{Moves on Gauss diagrams}
\label{fig:exrmoves}
\end{figure}

\begin{figure}[[htb] \epsfysize = 0.5 in
\centerline{\epsffile{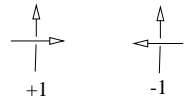}}
\caption{Crossing sign}
\label{fig:sign}
\end{figure}

An oriented  virtual  link diagram is a decorated immersion of $n$ oriented copies of $S^1 $ into the plane with two types of crossings: classical and virtual. The classical crossings are indicated by over/under markings and the virtual crossings are indicated by a circle around the crossings. 
The extended set of Reidmeister moves is shown in figure \ref{fig:exrmoves}. 
Two oriented virtual link diagrams are said to be equivalent if they are related by a finite sequence of the extended Reidmeister moves. An oriented virtual link is an equivalence class of  oriented virtual link diagrams. 
An oriented virtual knot is a single component virtual link. The sign of a classical crossing, $x$ is determined by the right hand rule (figure \ref{fig:sign}) and is denoted $sgn(x) $.

\begin{figure}[htb] \epsfysize = 1.0 in
\centerline{\epsffile{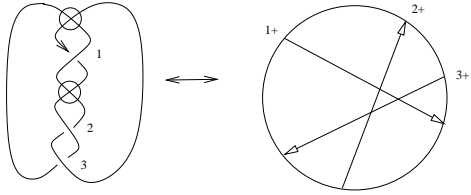}}
\caption{Corresponding diagrams}
\label{fig:knotequiv}
\end{figure}

There is a one to one correspondence between virtual knots and Gauss diagrams \cite{kvirt}.
\begin{thm}[Kauffman]\label{cor} Equivalence classes of Gauss diagrams are in one to one correspondence with oriented virtual knots. \end{thm}
We construct a Gauss diagram from an oriented virtual knot diagram by selecting a base point on the knot. We assign a label to each classical crossing.  As we traverse the knot, we record crossing information as a decorated symbol.  The symbol 
includes the label of the crossing traversed, a decoration indicating whether the overpass or underpass was traversed and the sign of the crossing. For a knot diagram with $n$ crossings, this results in a Gauss code with $2 n$ decorated symbols where each crossing label occurs twice. The decorated symbols are recorded in a clockwise orientation around a circle. The two occurrences of a label on the circle are connected by a signed oriented chord; the arrowhead placed at the overpassing symbol and the chord is marked with the sign of the crossing. We show a Gauss diagram and its corresponding knot (figure \ref{fig:knotequiv}).

\begin{figure}[htb] \epsfysize = 1.0 in
\centerline{\epsffile{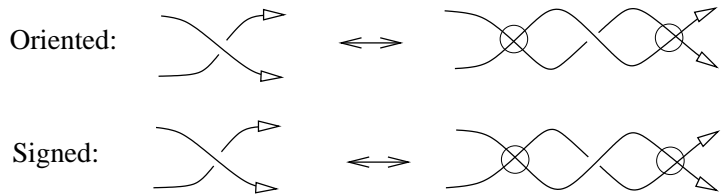}}
\caption{The diagrammatic virtualization moves}
\label{fig:virtualize}
\end{figure}

The chord moves are analogs of the Reidemeister moves.
The single chord move corresponds to a Reidemeister I move, the two chord move corresponds to the Reidemeister II move, and the two triangle move corresponds to the Reidemeister III move. There are no analogs of the moves involving virtual crossings since they do not alter the classical crossings. 
Crossing change in the knot diagram corresponds to changing the sign and orientation of the chord in the Gauss diagram. There are two methods of virtualization (see figure \ref{fig:virtualize}).  In oriented virtualization: the orientation of the chord is reverse in the Gauss diagram. In the corresponding knot diagram, the over passing strand is changed to the under passing strand in the crossing and the crossing is flanked by virtuals. The sign of the crossing does not change in orientated virtualization.  In signed virtualization, the sign of the chord is changed. In the corresponding knot diagram, we obtain an oppositely signed crossing flanked by virtuals. 
In figure \ref{fig:virtualize},  we see an example of virtualization. 

\begin{figure}[htb] \epsfysize = 2 in
\centerline{\epsffile{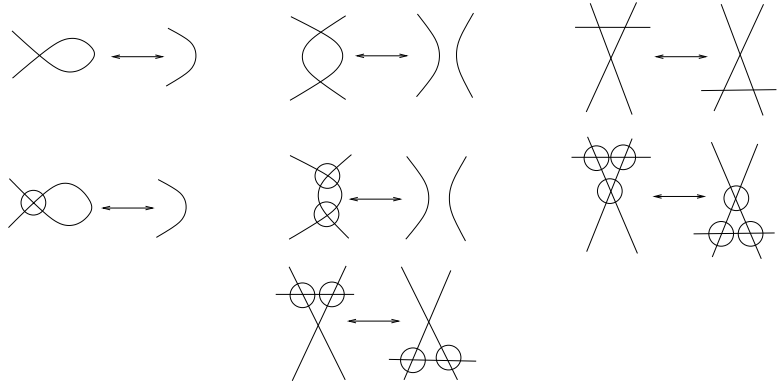}}
\caption{The flat diagrammatic moves}
\label{fig:flat}
\end{figure}

We define flat knots and links. A flat link diagram is a decorated immersion of $n$ copies of the circle into the plane with two types of crossings: flat and virtual. The flat crossings are indicated by a solid double point. Two flat virtual diagrams are said be equivalent if they are related by a sequence of the flat diagrammatic moves (figure \ref{fig:flat}). Flat links are equivalence classes of flat link diagrams. 
In the next section, we define parity mappings on Gauss diagrams and by extension virtual knot diagrams. 
\section{Parity mappings}
\begin{figure}[htb] \epsfysize = 1 in
\centerline{\epsffile{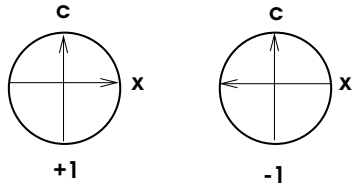}}
\caption{The value of $ int_c (x) $}
\label{fig:intx}
\end{figure}

Let $C$ denote the set of chords in a Gauss diagram, $G$.  Let $N_c$ be the set of chords that intersect chord $c$.  We define the intersection number of $x$ with $c$ (denoted $ int_c (x) $) as shown in figure
\ref{fig:intx}.

\noindent
A parity mapping on a Gauss diagram is a mapping $p: C \rightarrow \mathbb{Z} $ such that
\begin{equation} \label{paritydef}
p(c) = \sum_{x \in N_c } sgn(x) int_c (x) .
\end{equation}

\begin{lem}\label{chordinvar} The parity of existing chords is unchanged by the chord moves.
\end{lem}
\noindent
\textbf{Proof: } The one chord move introduces an isolated chord $h$. (The inverse one chord move removes an isolated chord.) The chord $h$ does not intersect any chords and has no effect on $ p(c)$ for any existing chords. 

The two chord move introduces a pair of oppositely signed chords $a$ and $b$. If a chord $c$ intersects $a$ then it also intersects $b$.  We observe that $ int_c (a) = int_c (b)$. Then:
\begin{equation*}
 sgn(a) int_c (a) + sgn(b) int_c (b) =0.
\end{equation*}
Hence, for any chord $c$, the value of $p(c)$ is unchanged. 

\begin{figure}[htb] \epsfysize = 1 in
\centerline{\epsffile{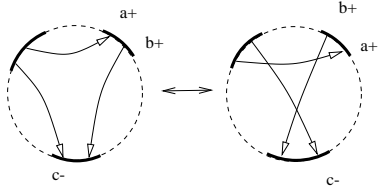}}
\caption{A (3,0) triangle move}
\label{fig:trianglemove1}
\end{figure}

We consider the triangle moves. No new chords are introduced during the triangle move. Chords and their intersections outside of the triangle of chords  are not affected. A (3,0) triangle move is 
shown in  \ref{fig:trianglemove1}.
We designate the chords as $a, b, $ and $c$. Let $ R_i $ denote the summands in $p(i)$ from chords other than
$a,  b,$ and $c$. We compute the intersection numbers for the zero side of the (3,0) move.
\begin{align*}
int_{a} (b) &= 0 &  int_{a}  (c) & =0 \\
int_{b} (a) & = 0 &  int_{b}  (c) &= 0\\
int_{c} (a) &=0 &  int_{c} (b) &=0
\end{align*}
The parity of the chords:
\begin{align*}
p(a) &= R_a,\\
p(b) &= R_b, \\
p(c) &= R_c.
\end{align*}

\begin{figure}[htb] \epsfysize = 1 in
\centerline{\epsffile{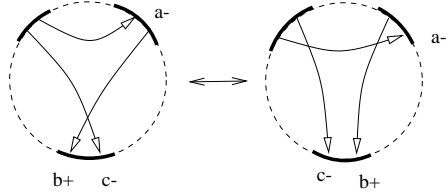}}
\caption{Triangle move 2}
\label{fig:trianglemove2}
\end{figure}

We compute the parity for the three intersection side of the (3,0) move. The intersection numbers have changed.
\begin{align*}
int_{a}(b)& = 1 &  int_{a}  (c) & =1 \\
int_{b} (a) & = 1 & int_{b} (c) & = 1\\
int_{c} (a) & =-1 &  int_{c} (b) &=1
\end{align*}
The parity of the chords $a, b,$ and $c$ are given by the formulas below.
\begin{align*} 
p(a) &= R_a + sgn(b) + sgn(c)\\
p(b) &= R_b + sgn(a) + sgn(c)  \\
p(c) &= R_c - sgn(a) + sgn(b)
\end{align*}
We note that $sgn(a) = sgn(b) = +1$ and $sgn(c)=-1$.
Hence, the parity of the chords is unchanged.
We analyze the (1,2) triangle move (figure \ref{fig:trianglemove2}).
We designate the chords as $a, b $ and $c$. 
We compute the intersection numbers on the left hand side of  figure \ref{fig:trianglemove2}.

\begin{align*}
int_{a} (b) &= 0 &  int_{a}  (c) &=0 \\
int_{b} (a) &= 0 & int_{b}  (c) &= -1\\
int_{c} (a) &=0 &  int_{c} (b) &=1
\end{align*}

The parity of the chords on the left hand side.
\begin{align*} 
p(a) &= R_a \\
p(b) &= R_b - sgn(c) \\
p(c) &= R_c + sgn(b)
\end{align*}

We compute the intersection number of the chords on the right hand side.
\begin{align*}
int_{a} (b)&= 1 & int_{a} (c)& =1 \\
int_{b} (a) &= -1 &  int_{b}  (c) &= 0\\
int_{c}  (a) &=-1 &  int_{c}  (b) & =0
\end{align*}

We compute the parity of the chords in the diagram on the right hand side.
\begin{align*} 
p(a) &= R_a + sgn(b) + sgn(c) \\
p(b) &= R_b - sgn(a)     \\
p(c) &= R_c - sgn(a)  
\end{align*}
Since $sgn(a) = sgn(c) = -1 $ and $sgn(b) = +1 $, we observe that the parity of the chords is unchanged by the diagrammatic move.
For the both the (3,0) and (1,2) triangle moves, there are three other oriented, signed pairs of  Guass diagrams that correspond to Reidemeister III moves. 
The first  is obtained by changing the sign and orientation of chord $a$. The second  is obtained by changing the sign and orientation of chords $b$ and $c$. The third is obtained by flipping the sign and orientation of chords $a$ and $c$. A short calculation shows that the triangle move does not change the parity of the chords in these diagrams. \qed

\begin{lem} \label{change} Let $G $ and $G'$ be Gauss diagrams with chord sets $C$ and $C'$ respectively. The set $C'$  is obtained  from $C$ by flipping the orientation and sign of the chord $y$  to form $y' $ in $C'$. 
We have generalized parity mappings:
$p: C \rightarrow \mathbb{Z} $ and $ p': C' \rightarrow  \mathbb{Z} $. 
\noindent
For $ c  \neq y, y' $:
\begin{equation*}
p(c) = p'(c).   
\end{equation*}
For $y $ and $y' $:
\begin{equation*}
p(y) = - p'(y').
\end{equation*} 
\end{lem}

\noindent
\textbf{Proof: }
For $c \neq y, y' $. If $c$ does not intersect $y$ in $G$ (respectively $y'$ in $G'$) 
then $ p(c) = p'(c)$. Suppose $c$ intersects $y$ in $G$ (respectively $y'$ in $G'$).
\begin{align*} 
p(c) &=int_c (y) sgn(y)  +  \sum_{x \in N_c , x \neq y}  int_c (x) sgn (x) \\
p'(c) &=int_c (y') sgn(y')   + \sum_{x \in N_{c}, x \neq y' } int_{c} (x) sgn(x)
\end{align*}
Since $ int_c (y) sgn(y) = int_c (y') sgn(y')$ this implies that $p(c) = p'(c) $. 
We consider the parity of $y$ and $y'$.
\begin{align*} 
p(y) &=  \sum_{x \in N_y}  int_y (x) sgn (x) \\
p'(y') &=  \sum_{x \in N_{y'}} int_{y'} (x) sgn(x)
\end{align*}
For $y$ and $y'$, we note that $N_y = N_{y'} $ and  that $int_y (x) =- int_{y'} (x) $. 
Hence, $p(y) = - p'(y') $. 
\qed

\begin{rem}
From the generalized parity mapping, we obtain an element of  $ \mathbb{Z} $ for each chord. Hence, we associate an element to each crossing in the corresponding knot diagram. We can compose the parity mapping with the projection $ \phi_n : \mathbb{Z} \rightarrow \mathbb{Z}_n $ to associate an element of $ \mathbb{Z}_n $ with each chord. If we let $n=2$, we obtain Manturov's original definition of parity \cite{m2}. 
We can obtain mappings into non-cyclic group by assigning generators to each chord in such a way that the assignment respects the analogs of the Reidemeister moves. This leads to non-trivial mappings for classical knots. 
\end{rem}

\section{Gauss diagram invariants}
Let $i$ be an element of $ \mathbb{Z} $. For a gauss diagram $G$, let $C$ be the set of chords in
$G$. We define
 $A_i  (G) $ as
\begin{equation}
A_i (G) = \lbrace c \in C | p(c) = i \rbrace.
\end{equation}

The set $A_i (G) $  is not an invariant of the equivalence class of the Gauss diagram $G$.
We define 
the signed cardinality of the set $A_i (G) $ (denoted $|A_i (G) | $) to be:
\begin{equation} 
|A_i (G)| = \sum_{c \in A_i (G) } sgn(c) .
\end{equation}
The constants $ | A_i (G) | $, $i \in \mathbb{Z}^* $ are invariants of equivalence classes of Gauss diagrams.

\begin{thm}For a Gauss diagram $G$ and a non-zero integer $i$,  $ | A_i (G) |$ is invariant under 
the single chord move, two chord move and the triangle move. \end{thm}
\noindent
\textbf{Proof:} The single chord move introduces a chord with zero parity. The two chord move introduces two chords with the same parity. However, these chords are oppositely signed and make a net contribution of zero to the signed cardinality. The triangle move does not introduce any new chords and the parity of the existing chords is unchanged by the triangle move. \qed

We now consider n-tuples. Let $ Z  = (z_1, z_2, \ldots z_n) $ such
 that $z_i \neq 0 $ and $z_i  < z_{i+1} $. We define a set of n-tuples of chords
\begin{equation}
 A_{Z}  (G)  = \lbrace (x_1, x_2, \ldots x_n) |  x_i \in A_{z_i}  (G) \text{ for }  i \in 1,2, \ldots n 
\rbrace .
\end{equation}

Let $ \bar{X}  = (x_1, x_2, \ldots x_n ) $ denote an n-tuple of chords. The signed cardinality 
of this set  is 
\begin{equation}
|A_{Z}  (G)  | = \sum_{\bar{X} \in A_{ Z}} sgn(x_1) sgn(x_2) \ldots sgn(x_n).
\end{equation}

\begin{thm} \label{first} The sum $ |A_{Z} (G) | $ is an invariant of the Gauss diagram. \end{thm} 

\noindent
\textbf{Proof:} The signed cardinality is invariant under the single chord move, since $A_0 (G) $ does not contribute any chords. The two chord move does not change $ | A_Z (G) | $ since the contributions from the oppositely signed pair of chords in the two chord move cancel each other.  Parity does not change under the triangle move. \qed

Let $i$ be a non-negative integer. We define
\begin{equation}  \label{vset1}
V_i (G) = \lbrace x \in C | |p(x)|=i \rbrace .
\end{equation}
 We define $|V_Z (G) |$
\begin{equation} \label{sgnvset1}
|V_i (G)|  = \sum_{x \in V_{i} (G) } sgn(x).
\end{equation}
Now, let $Z = (z_1, z_2, \ldots z_n) $ such that $ 0< z_1 $ and $ z_i < z_{i+1} $.
\begin{equation} \label{vset}
 V_{Z}  (G)  = \lbrace (x_1, x_2, \ldots x_n) |  x_i \in V_{z_i} (G)  \text{ for }  i \in 1,2, \ldots n 
\rbrace .
\end{equation} \label{sgnvset}
 We define $|V_Z (G) |$ using equation \ref{vset}
\begin{equation}
|V_Z (G)| =  = \sum_{\bar{X} \in V_{ Z} (G) } sgn(x_1) sgn(x_2) \ldots sgn(x_n).
\end{equation}
\begin{thm} \label{second} The sums $|V_i (G)| $ (equation \ref{sgnvset1}) and $ |V_Z (G) | $ (equation \ref{sgnvset}) are invariant under the chord diagram moves. \end{thm} 

\noindent
\textbf{Proof:} The sum $ |V_i (G) | = |A_i (G)| + |A_{-i} (G)| $.  The $ |V_Z (G) | $ can be expressed as a linear combination of $ |A_{Z} (G) |$ .\qed

\section{Vassiliev Invariants of Knots}

Equivalence classes of virtual knot diagrams and equivalence classes of Gauss diagrams are in a one to one correspondence. For any oriented virtual knot diagram $K$, there is a corresponding Gauss diagram $G_K$.  Each classical crossing in $K$ corresponds to an oriented, signed chord in $G_K$. We partition the set of crossings in
$K$ using this correspondence. 
If a crossing $c$ corresponds to a chord with parity $i$ then we say that the crossing has parity $i$.

We define $ A_i (K) $ to be the set of crossings in $K$ that correspond to chords with parity $i$ in 
$A_i (G_K) $.  We analogously define $V_i (K) $ to be the set of crossings with parity $i$ or $-i$. Hence,   the n-tuples of crossings in the sets: $A_{Z} (K) $ (respectively $ V_{ Z } (K) $) correspond to the  n-tuples of chords in the sets: $ A_{Z} (G_K) $ (respectively $V_{ Z} (G_K) $). 
We immediately obtain the following theorem.
\begin{thm} The sums $ |A_i (K)| $, $| A_Z (K)| $, $ |V_i (K)| $, and $| V_Z (K)| $ are invariant under the extended set of Reidemeister moves. \end{thm}

We define a formal sum of flat knot diagrams following the methods from Henrich \cite{allison}.
However, we first need to introduce notation that will allow us to indicate singular crossings, resolutions of singular crossings, vertically smoothed crossings, and flattened diagrams.
A singular crossing in a virtual knot is a rigid vertex that the extended Reidemeister moves do not apply to. 
 Let $\bar{x} = (x_1, x_2, \ldots x_n)  $ denote $n$ crossings in the knot $K$. The notation $K_x $ denotes a knot with a one singular crossing,  $K_{\bar{x} } $ denotes a 
knot with singular crossings at the crossings: $ x_1, x_2, \ldots x_n$.  

Let $ \bar{c}  = (c_1, c_2, \ldots c_n ) $ denote an element of $ \lbrace \pm 1 \rbrace ^n$. The notation 
$K_{ (\bar{c}, \bar{x} ) } $ indicates the resolution of the $n$ singular crossings in $K_{\bar{x}} $. The
singular crossing $x_i $ has been resolved as a positive crossing (respectively negative) if $ c_i = +1 $
($c_i = -1 $).  For the case with one or two singular crossings, we write $ K_{(\pm 1,x)} $ or $K_{( \pm \pm, xy)} $.

The notation $F(K) $ indicates a flat knot obtained from $K$ by  flattening the classical crossings in $K$.
The notation $K^z $ or $ K^{\bar{z}} $ indicates that the individual crossing $z$ or the n-tuple of 
crossings $\bar{z} $ have been smoothed vertically.

We define $sgn(\bar{c} )$.
\begin{equation*}
sgn(\bar{c} ) = \prod_{i=1} ^{n} c_i .
\end{equation*}
For the $n$ classical crossings $ \bar{z} $, we define 
\begin{equation*}
sgn(\bar{z}) = \prod_{i=1} ^{n} sgn(z_i) 
\end{equation*}

We recall the definition of a Vassiliev invariant. Aknot invariant $ \mathcal{S} $ is a Vassiliev invariant of degree $n$ if $ \mathcal{S}(K) = 0 $ for any knot $K$ with $n+1 $ or more singular crossings \cite{gpv}. From this definition, we obtain
\begin{equation}\label{step1}
\mathcal{S} (K_{x} ) = \mathcal{S} (K_{(+,x)} ) - \mathcal{S} (K_{(-,x)}) 
\end{equation}

For a knot with $n$ singluarities, from equation \ref{step1}, we obtain
\begin{equation} \label{stepn}
\mathcal{S} (K_{\bar{x}}) = \sum_{\bar{c}  \in \lbrace \pm 1 \rbrace ^n } sgn( \bar{c})  \mathcal{S}  (K_{ (\bar{c} ,\bar{x} )}).
\end{equation} 

We define two families of  Vassiliev invariants below and use equations \ref{step1} and \ref{stepn} to verify the
degree of these invariants.  

We define a family of degree one Vassiliev invariants using the parity of the crossings. 
\begin{equation}
S_{i} (K) = \sum_{ x \in V_{i} (K)  } sgn(x) F(K^x )
 \end{equation}
where $i \in \mathbb{N}$. 
This invariant is a formal sum of flat diagrams. The summands are obtained by smoothing a single crossing with parity $i$ and assigned a weight based on the sign of that crossing.
This results in formal sums of flat diagrams with coefficients in $ \mathbb{Z} $.

\begin{thm}For $i \neq 0 $,  $S_i (K) $ is a degree one Vassiliev invariant. \end{thm}

\noindent
\textbf{Proof:}
The formal sum is unchanged by the extended set of Reidemeister moves.
Recall that the parity of the existing crossings is unchanged by the Reidemeister moves. 
Crossings introduced by a Reidemeister I move have parity zero and do not contribute to 
this formal sum. If $x$ and $y$ are two crossings involved in a Reidemeister II move, $ F(K_x) $ is equivalent as a flat to $F( K_y) $ and the net contribution is to the formal sum is zero.  
The Reidemeister III move does not change parity and the flat diagrams obtained from both sides are equivalent. 

We claim that 
$ S_i (K) $ is a degree one Vassiliev invariant.  Let $K_{ab} $ be a knot with singularities at $a, b $ .
We apply the definition of a Kaufman finite type invariant given in equation \ref{step1}  to show that $S_i $ vanishes on any knot with two singularities.  
\begin{equation} \label{write1}
S_{i} (K_{ab})  = S_{i} (K_{(++,ab) } - K_{(+-,ab)} - K_{(-+,ab)} + K_{(--,ab)} )
\end{equation}
 By Lemma \ref{change}, if $a+ \in  V_i (K) $ then $ a- \in V_i (K)$. We assume that $a,b  \in V_i (K) $.

Expanding equation \ref{write1}:
\begin{multline*}
S_{i} (K_{ab}) = F(K_{(++,ab)} ^a ) + F(K_{(++,ab)} ^b) \\
 - [ F(K_{(+-,ab)}^a) - F(K_{(+-,ab)}^b)] -
 [-F(K_{-+,ab)}^a) + F(K_{(-+,ab)}^b )] \\
 + [ -F(K_{(--,ab)}^a) -F (K_{(--,ab)}^b) ]
=0 .
\end{multline*}
If $a$ is an element of $V_i (K) $ and $b$ is not an element of $V_i (K) $ then 
\begin{multline*}
S_{i} (K_{ab}) = F(K_{(++,ab)} ^a )  - F(K_{(+-,ab)} ^a) \\
+ F(K_{(-+,ab)}^a)  -F(K_{(--,ab)}^a)
=0 .
\end{multline*}
If $a$ and $b$ are not elements of $V_i (K)$ then $ S_i (K^{ab}) =0 $. \qed
\begin{rem} Henrich's smoothed, degree one Vassiliev invariant \cite{allison} can be decomposed as 
a sum of the the $S_i (K)$. 
\end{rem}
We define a family of degree $n$ Vassiliev invariants using the parity of the crossings. 
Let $Z$ be an element of $ \mathbb{N}^n$ with $ z_i < z_{i+1} $. 
\begin{equation}
S_{Z} (K) = \sum_{ \bar{x} \in V_{Z} (K)  } sgn(\bar{x}) F(K^{ \bar{x}}). 
 \end{equation}
The degree $n$ invariants are a formal sum of flat diagrams with coefficients in $ \mathbb{Z} $.
\begin{thm}$ S_{Z} (K) $ is a degree $n$ Vassiliev invariant. \end{thm}

\noindent
\textbf{Proof:} The formal sum of flat diagrams is invariant under the Reidemeister moves. Let  $ \bar{x} $ be an n-tuple of crossings. From equation  \ref{stepn}:
\begin{equation} \label{sings}
S_{Z} (K_{\bar{x}}) = \sum_{\bar{c} \in \{ \pm 1 \} ^n } sgn (\bar{c}) S_{Z} (K_{(\bar{c}, \bar{x} )} ).
\end{equation} 
For $ S_{Z} ( K_{ ( \bar{c}, \bar{x} )}  ) $ from equation \ref{sings}:
\begin{equation} \label{individ}
S_{Z} (K_{ ( \bar{c}, \bar{x} ) } )  = \sum_{\bar{z} \in V_{Z} (K_{(\bar{c}, \bar{x})}) } sgn(\bar{z} )
 K_{ (\bar{c}, \bar{x} )} ^{\bar{z}}.
\end{equation}
Hence, from equations \ref{sings} and \ref{individ}
\begin{equation} \label{oneway}
S_{Z} (K_{  \bar{x}  } )  = \sum_{\bar{c} }  \sum_{ \bar{z} \in V_{Z} (K_{(\bar{c}, \bar{x})} )} sgn (\bar{c}) sgn(\bar{z} ) 
 K_{ (\bar{c}, \bar{x} )} ^{\bar{z}}.
\end{equation}
Applying Lemma \ref{change} and rewriting equation \ref{oneway}
\begin{equation}
S_Z (K_{\bar{x}}) = \sum_{ \bar{z} \in V_{Z} (K_{ \bar{x}})            }  \sum_{\bar{c} }   sgn (\bar{c}) sgn(\bar{z} ) 
 K_{ (\bar{c}, \bar{x} )} ^{\bar{z}}.
\end{equation}

If $ \bar{z} $ is an element of $V_Z (K_{ \bar{x} } ) $, then we have three possibilities. 
\begin{enumerate}
\item There is no overlap 
between the crossings in $ \bar{z} $ and $ \bar{x}$. 
\item The n-tuples $ \bar{z} $ and $  \bar{x} $ have $i$ crossings in common where $0< i < n $. 
\item The n-tuples $ \bar{z} $ and $ \bar{x} $ are equivalent.
\end{enumerate}

\noindent
\textbf{Case 1:}

If there is no overlap between $ \bar{z} $ and $ \bar{x} $ then the net contribution to $ S_Z (K_{\bar{x}}) $ is zero.

\noindent
\textbf{Case 2:}

If the n-tuples $ \bar{z} $ and $ \bar{x} $ have crossings in common, then we decompose $ sgn(\bar{z})  $ as $ sgn(\bar{z} ') sgn( \bar{c}')  $.  The notation $sgn(\bar{z}') $ indicates the product of the signs of the crossings that are not in $ \bar{c} $. The notation $sgn(\bar{c}')$ indicates the product of the signs of the crossings in $ \bar{z} $ that are also in $ \bar{c} $. 

 The contribution to $S_Z (K_{ ( \bar{c}, \bar{x} )}  ) $ is
\begin{gather*}
\sum_{\bar{c} \in \lbrace \pm 1 \rbrace ^n } sgn(\bar{c})  sgn(\bar{z})  F( K_{ ( \bar{c}, \bar{x} )} ^{ \bar{z}} )  \\
= sgn(\bar{z}') \sum_{\bar{c} \in \lbrace \pm 1 \rbrace ^n} sgn(\bar{c})  sgn(\bar{c}')  F( K_{ ( \bar{c}, \bar{x} )} ^{ \bar{z}} )  \\
\end{gather*}
The term $ sgn( \bar{c} ) sgn( \bar{c}' ) $ simply squares some of the $c_i $. Hence, the net contribution is zero.

\noindent
\textbf{Case 3:}

If the n-tuples $\bar{z} $ and $ \bar{x}$ are equivalent, then the contribution to  $S_Z (K_{\bar{x}} ) $ is
\begin{gather*}
\sum_{\bar{c} \in \lbrace \pm 1 \rbrace ^n } sgn(\bar{c})  sgn(\bar{x})  F( K_{ ( \bar{c}, \bar{x} )} ^{ \bar{x}} )  \\
= \sum_{\bar{c} \in \lbrace \pm 1 \rbrace ^n} sgn(\bar{c})  sgn(\bar{c})  F( K_{ ( \bar{c}, \bar{x} )} ^{ \bar{x}} )  \\
= \sum_{ \bar{c} \in \lbrace \pm 1 \rbrace ^n}   F( K_{ ( \bar{c}, \bar{x} )} ^{ \bar{x}} ) .
\end{gather*}
The net contribution is $2^n $ copies of the diagram $ F (K ^{ \bar{x}} ) $.
We conclude that if $ \bar{x} \in V_Z (K) $ then
\begin{equation*}
S_Z (K_{\bar{x}}) = 2^n F(K^{ \bar{x}} ).
\end{equation*}
Hence, 
\begin{equation*}
S_Z (K_{a \bar{x}}) = 2^n F(K _{(+,a) }^{\bar{x} }) - 2^n F(K_{(-,a) }^{\bar{x} }) = 0
\end{equation*}
As a result, $S_Z (K) $ vanishes on any knot with $n+1 $ or more singularities.\qed

\begin{cor} The invariant $ |V_Z (K) | $ is a Vassiliev invariant. \end{cor}

\noindent 
\textbf{Proof:} We obtain $ |V_Z (K) | $ from $ S_Z (K) $  by viewing each flat diagram as a state and evaluating each state as $1$. \qed

\section{Examples}
\subsection{Virtual Trefoil}
\begin{figure}[htb] \epsfysize = 1.5 in
\centerline{\epsffile{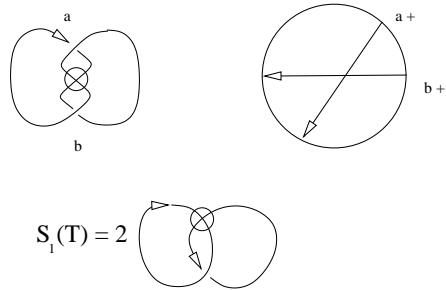}}
\caption{Virtual trefoil and corresponding Gauss diagram}
\label{fig:virtualtrefexa}
\end{figure}

We consider the right handed virtual trefoil, $T$(figure \ref{fig:virtualtrefexa}). Both crossings have positive sign. 
We compute the parity of the crossings from the Gauss diagram.
\begin{gather*} 
p(a) = 1  \text{ and } p(b) =1 
\end{gather*}
The value of $ S_1 (T) $ is displayed in figure \ref{fig:virtualtrefexa} and $ |V_1 (T)| = 2 $.

\subsection{Kishino's knot}

\begin{figure}[htb] \epsfysize = 1 in
\centerline{\epsffile{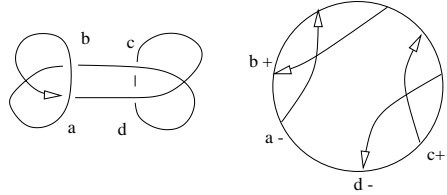}}
\caption{Kishino's knot}
\label{fig:kishino}
\end{figure}
For Kishino's knot, $K$: $ p(a) = p(b) = -1 $ and $p(c) = p(d) = 1 $. All crossings contribute to the formal sum, and $ |V_1 (K)| =0$. 

\subsection{Miyazawa knot}
\begin{figure}[htb] \epsfysize = 1.5 in
\centerline{\epsffile{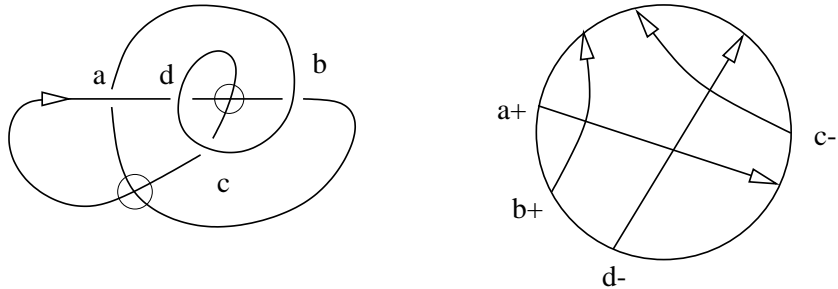}}
\caption{Miyazawa's knot}
\label{fig:miyazawa}
\end{figure}

\begin{figure}[htb] \epsfysize = 2.2 in
\centerline{\epsffile{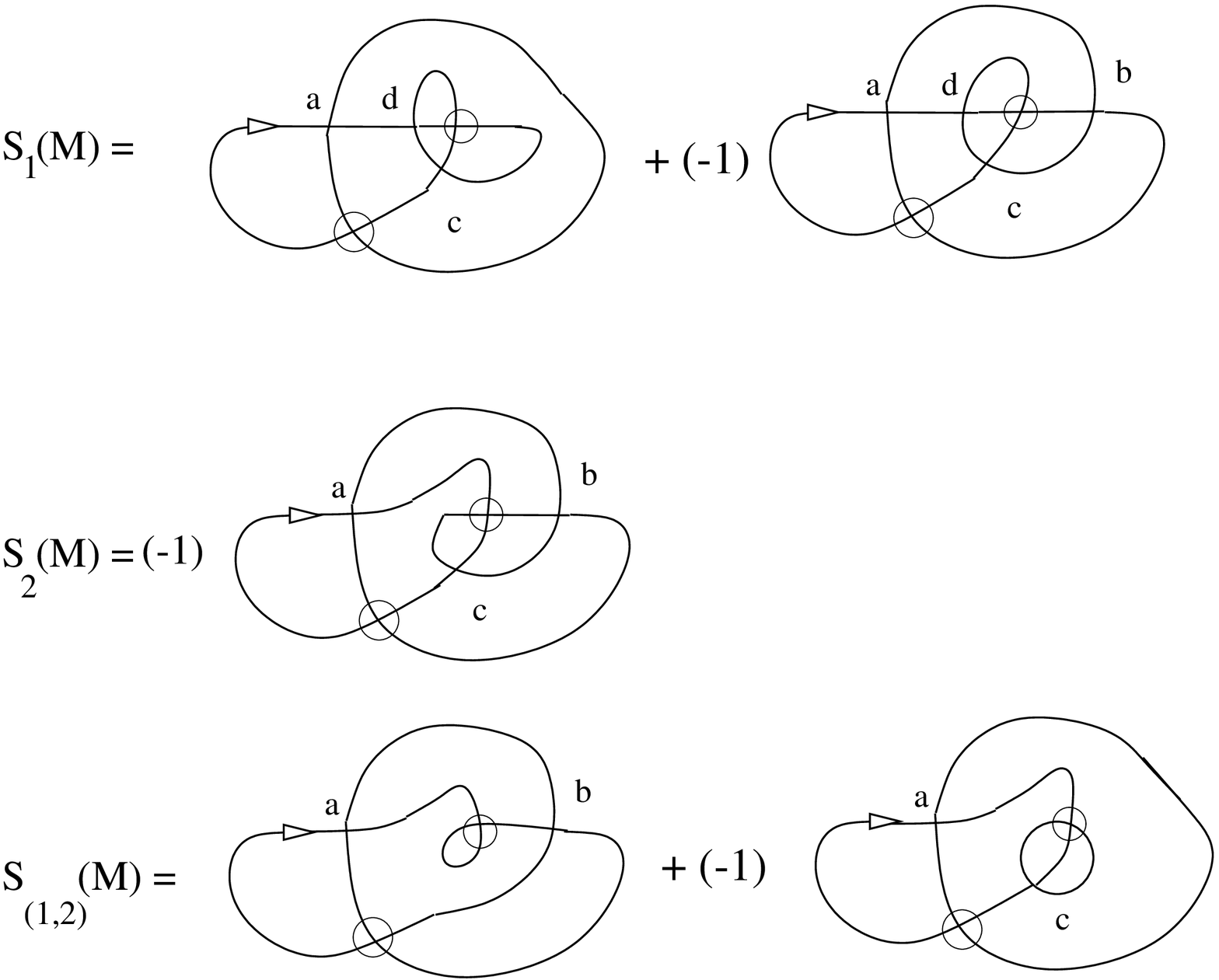}}
\caption{$ S_Z (M) $}
\label{fig:miyazawaevals}
\end{figure}

Let $M$ denote the Miyazawa knot (figure \ref{fig:miyazawa}).
We observe that $p(a) =0$, $p(b) = 1 $, $p(c) = -1 $ and $p(c) = 2$. 
For this knot, $ S_1 (M) $ and $S_2 (M) $ are non-zero. We show these formal sums in 
figure \ref{fig:miyazawaevals}. 
The formal sum $ S_{(1,2) } $ has two summands. The summands are the links obtained by smoothing the pairs $(c,d) $ and $(b,d) $ with the appropriate sign. The formal sum is shown in 
figure \ref{fig:miyazawaevals}. 
\subsection{Pretzel knot}
\begin{figure}[htb] \epsfysize = 1.5 in
\centerline{\epsffile{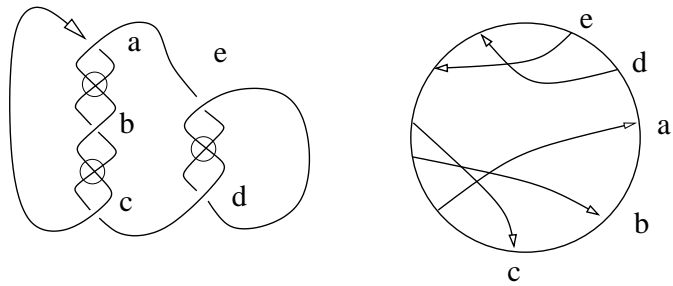}}
\caption{Pretzel knot}
\label{fig:pretzel}
\end{figure}

Let $P$ denote the pretzel knot (figure \ref{fig:pretzel}).
\begin{figure}[htb] \epsfysize = 2.2 in
\centerline{\epsffile{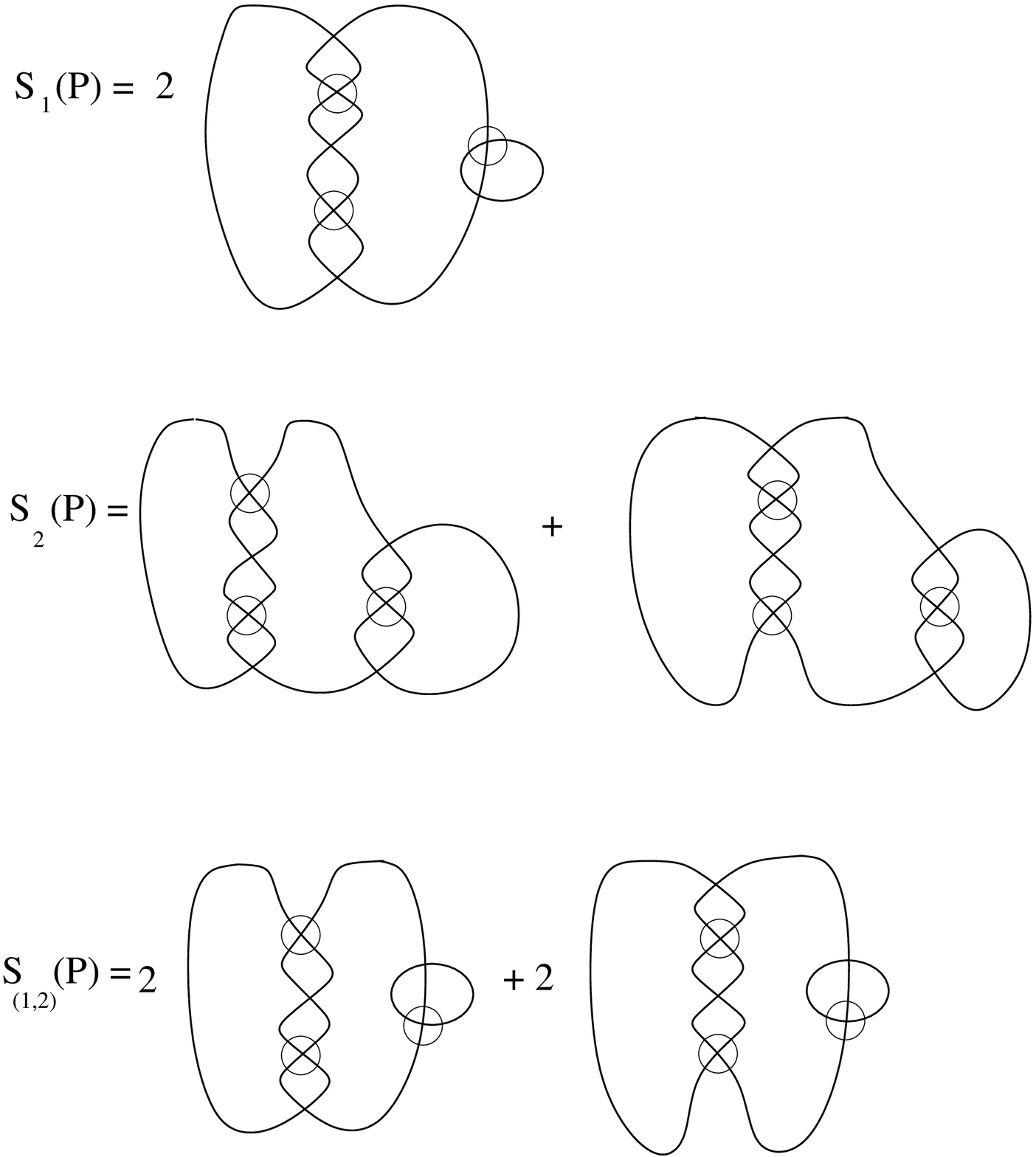}}
\caption{$S_Z (P) $}
\label{fig:pretzelevals}
\end{figure}

We observe that $p(a) =2$, $p(b) = 0 $, $p(c) = -2 $ , $p(d) = 1$, and $p(e) = -1 $. 
For this knot, $ S_1 (M) $ and $S_2 (M) $ are non-zero. We show these formal sums in 
figure \ref{fig:pretzelevals}. 
The formal sum $ S_{(1,2) } $ has two summands.  The formal sum is shown in 
figure \ref{fig:pretzelevals}.

\end{document}